\documentclass[12pt]{article}
\usepackage{graphicx}
\usepackage{amsmath}
\usepackage{graphics}
\usepackage{epsfig}
\usepackage{amssymb}

\newcommand{\be}{\begin{eqnarray}}
\newcommand{\ee}{\end{eqnarray}}
\newcommand{\bea}{\begin{eqnarray*}}
\newcommand{\eea}{\end{eqnarray*}}

\def\ce{\mathbb{C}}
\renewcommand{\epsilon}{\varepsilon}
\def\eps{\varepsilon}
\def\er{\mathbb{R}}


\begin{document}
\title{Random block matrices and matrix orthogonal polynomials}

\author{
{\small Holger Dette} \\
{\small Ruhr-Universit\"at Bochum} \\
{\small Fakult\"at f\"ur Mathematik} \\
{\small 44780 Bochum, Germany} \\
{\small e-mail: holger.dette@rub.de}\\
{\small FAX: +49 234 3214 559}\\
\and
{\small Bettina Reuther} \\
{\small Ruhr-Universit\"at Bochum} \\
{\small Fakult\"at f\"ur Mathematik} \\
{\small 44780 Bochum, Germany} \\
{\small e-mail: bettina.reuther@rub.de}\\
{\small FAX: +49 234 3214 559}\\
}

\maketitle

\begin{abstract}
In this paper we consider random block matrices, which generalize the general beta ensembles, which were recently investigated by Dumitriu and
Edelmann (2002, 2005). We demonstrate that the eigenvalues of these random matrices can be uniformly approximated by roots of matrix orthogonal
polynomials which were investigated independently from the random matrix literature. As a consequence we derive the asymptotic spectral
distribution of these matrices. The limit distribution has a density, which can be represented as the trace of an integral of densities of
matrix measures corresponding to the Chebyshev matrix polynomials of the first kind. Our results establish a new relation between the theory of
random block matrices and the field of matrix orthogonal polynomials, which have not been explored  so far in the literature.

\end{abstract}


\section{Introduction}
\def\theequation{1.\arabic{equation}}
\setcounter{equation}{0}

The classical orthogonal, unitary and symplectic ensembles of random matrices have been studied extensively in the literature, and there are
numerous directions to generalize these investigations to the case of random block matrices. Several authors have studied the limiting spectra
of such matrices under different conditions on the blocks
 [see e.g.\
Girko (2000), Oraby (2006 a,b) among others].

In the present paper we study a class of random block matrices which
generalize the tridiagonal matrices corresponding to the classical
$\beta$-Hermite ensemble defined by its density
\begin{equation}\label{1.1}
c_\beta \cdot \prod_{1 \le i < j \le n} \vert \lambda_i - \lambda_j\vert^\beta e^{- \sum^n_{i=1}\frac{\lambda^2_i}{2}} ~,
\end{equation}
where $\beta > 0$ is a given parameter and $c_\beta$ a normalizing constant. It is known for a long time that for $\beta = 1,2,4$ this density
corresponds to the Gaussian ensembles [see Dyson (1962)], which have been studied extensively in the literature on random matrices [see Mehta
(1967)]. It was recently shown by Dumitriu and Edelmann (2002) that for any $\beta > 0 $ the joint density of the eigenvalues $\lambda_1 \leq
\dots \leq \lambda_n$ of the $n \times n$ symmetric matrix
\begin{equation}\label{1.2}
G_n^{(1)}=\left(
\begin{matrix}
N_1 & \frac 1{\sqrt 2} X_{(n-1)\beta}\\
\frac 1{\sqrt 2}X_{(n-1)\beta} & N_2 & \frac 1{\sqrt 2}X_{(n-2)\beta}\\
 & \ddots & \ddots & \ddots \\
& & \frac 1{\sqrt 2} X_{2\beta} & N_{n-1} & \frac 1{\sqrt 2}X_\beta\\
& & & \frac 1{\sqrt 2}X_\beta  & N_n
\end{matrix}\right)
\end{equation}
is given by (\ref{1.1}) where $X_\beta, \dots, X_{(n-1)\beta}$, $N_1, \dots, N_n$ are independent random variables with $X^2_{j \beta} \sim
\mathcal{X}^2_{j \beta}$ and $N_j \sim \mathcal{N} (0,1)$. Here $\mathcal{X}^2_{j \beta}$ and $\mathcal{N} (0,1)$ denote a chi-square  (with $j
\beta$ degrees of freedom) and a standard normal distributed random variable, respectively. In the present paper we consider
 random block matrices of the form
\begin{eqnarray}\label{1.3}
G_n^{(p)}=\left(
\begin{array}{cccccccc}
B_{0}^{(p)} &  A_{1}^{(p)} &  &    &          \\
A_{1}^{(p)} & B_{1}^{(p)} &  A_{2}^{(p)} &   &   \\
 &\ddots  & \ddots & \ddots & &  \\
  &  &  A_{\frac{n}{p}-2}^{(p)}  &  B_{\frac{n}{p}-2}^{(p)} &  A_{\frac{n}{p}-1}^{(p)}  \\
 & & & A_{\frac{n}{p}-1}^{(p)} & B_{\frac{n}{p}-1}^{(p)} \\
 \end{array}
\right),
\end{eqnarray}
where $n = mp \in \mathbb{N}, ~m \in \mathbb{N}, ~p \in
\mathbb{N}$ and the $p \times p$ blocks $A^{(p)}_i$ and $ B^{(p)}_i$
are defined by
 \begin{eqnarray}\label{1.4}
B_i^{(p)}=\frac{1}{\sqrt{2}}\left(
\begin{array}{cccccccc}
\sqrt{2} N_{ip+1} &  X_{\gamma_1 (n-ip-1)} &X_{\gamma_2 (n-ip-2)}  & \cdots   &   X_{\gamma_{p-1} (n-(i+1)p+1)}       \\
X_{\gamma_1 (n-ip-1)} &\sqrt{2}N_{ip+2}& X_{\gamma_1 (n-ip-2)} & \cdots &    X_{\gamma_{p-2} (n-(i+1)p+1)}  \\
  &  &  &  &   \\
 \vdots &\ddots & \ddots & \ddots & \vdots  \\
  & & &  &   \\
X_{\gamma_{p-2} (n-(i+1)p+2)}  &\cdots  &  X_{\gamma_{1} (n-(i+1)p+1)} & \sqrt{2}N_{(i+1)p-1} & X_{\gamma_1 (n-(i+1)p+1)} \\
X_{\gamma_{p-1} (n-(i+1)p+1)}  & \cdots&  X_{\gamma_{2} (n-(i+1)p+1)} & X_{\gamma_1 (n-(i+1)p+1)} & \sqrt{2}N_{(i+1)p} \\
\end{array}
\right)
\end{eqnarray}
and
\begin{eqnarray}\label{1.5}
A_i^{(p)}=\frac{1}{\sqrt{2}}\left(
\begin{array}{cccccccc}
X_{\gamma_p(n-ip)} &  X_{\gamma_{p-1}(n-ip)}  & X_{\gamma_{p-2}(n-ip)}  & \cdots   &  X_{\gamma_{1}(n-ip)} \\
X_{\gamma_{p-1}(n-ip)} &X_{\gamma_p(n-ip-1)}& X_{\gamma_{p-1}(n-ip-1)}  & \cdots &   X_{\gamma_{2}(n-ip-1)}\\
  &  &  &  &   \\
 \vdots &\ddots & \ddots & \ddots & \vdots  \\
  & & &  &   \\
 X_{\gamma_2(n-ip)}  &\cdots  &  X_{\gamma_{p-1} (n-(i+1)p+3)}   & X_{\gamma_p (n-(i+1)p+2)} & X_{\gamma_{p-1} (n-(i+1)p+2)} \\
X_{\gamma_1(n-ip)}  & \cdots&  X_{\gamma_{p-2} (n-(i+1)p+3)}& X_{\gamma_{p-1} (n-(i+1)p+2)} & X_{\gamma_p (n-(i+1)p+1)} \\
 \end{array}
\right),
\end{eqnarray}
respectively,  $\gamma_1, \dots, \gamma_p > 0$ are given constants and all random variables are independent. Note that in the case $p=1$ we
have $G_n^{(p)} = G_n^{(1)}$ and that for $p>1$ the elements $g_{ij}^{(p)}$ of the matrix $G_n^{(p)}$ are standard normally distributed if
$i=j$ and $\mathcal{X}_k$-distributed else, where the degrees of freedom depend on the position of the element in the matrix $G_n^{(p)}$.
 In the present
paper we investigate the limiting spectral behaviour of matrices of the form (\ref{1.3}). In particular it is demonstrated that the eigenvalues
of the matrix $G_n^{(p)}$ are closely related to roots of matrix orthogonal polynomials, which have been studied independently from the theory
of random matrices [see e.g.\ Duran (1996, 1999), Duran and Daneri-Vias (2001), Duran and Lopez-Rodriguez (1996), Duran, Lopez-Rodriguez and
Saff (1999), Sinap and van Assche (1996) or Zygmunt (2002) among others]. To our knowledge this is the first paper, which relates  the field of
matrix orthogonal polynomials to the theory of random matrices.

In Section 2 we review some basic facts on matrix orthogonal polynomials. We prove a new result on the asymptotic behaviour of the roots of
such polynomials which is of own interest and   the basis for the investigation of the limiting spectral properties of the matrix $G^{(p)}_n$.
  In Section~3 we provide a strong uniform
approximation of the random eigenvalues of the matrix $G_n^{(p)}$ by the deterministic roots of matrix valued orthogonal polynomials, and these
results are applied  to study the asymptotic behaviour  of the spectrum of the matrix $G_ n^{(p)}$. In particular, we derive a new class of
limiting spectral distributions which generalize the classical Wigner semi circle law in the one-dimensional case. Roughly speaking, the limit
distribution has a density, which can be represented as the trace of an integral of densities of matrix measures corresponding to the Chebyshev
matrix polynomials of the first kind. Finally some examples are presented in Section 4 which illustrate the theoretical results.

\section{Matrix orthogonal polynomials}
\def\theequation{2.\arabic{equation}}
\setcounter{equation}{0}

In the following discussion we consider for each $k \in \mathbb{N}$ a sequence of $p \times p$ matrix polynomials $(R_{n,k}
(x))_{n \geq 0}$, which are defined recursively by
\begin{eqnarray}\label{2.1}
x R_{n,k}(x)=A_{n+1,k}R_{n+1,k}(x)+B_{n,k}R_{n,k}(x)+A^T_{n,k}R_{n-1,k}(x),
\end{eqnarray}
with initial conditions $R_{-1,k}(x)=0,~R_{0,k}(x)=I_p$, where $B_{i,k} \in \mathbb{R}^{p \times p}$ denote symmetric and $A_{i,k} \in
\mathbb{R}^{p \times p}$ denote non-singular matrices. Here and throughout this paper $I_p$ denotes the $p \times p$ identity matrix and 0 a $p
\times p$ matrix with all entries equal to 0. A matrix measure $\Sigma$ is a  $p \times p$ matrix of  (real) Borel measures such that  for each
Borel set $A \subset \mathbb{R}$ the matrix $\Sigma (A) \in \mathbb{R}^{p \times p}$ is nonnegative definite. It follows from Sinap and van
Assche (1996) that there exists a positive definite matrix measure $\Sigma_k \in \mathbb{R}^{p \times p}$   such that the polynomials $\left(
R_{n,k} (x)\right)_{n \geq 0}$ are orthonormal with respect to the inner product induced by the measure $d \Sigma_k (x)$, that is
\begin{eqnarray}\label{2.2}
\int R_{n,k}(x)d\Sigma_k(x)R_{m,k}^T(x)=\delta_{nm}I_p.
\end{eqnarray}
The roots of the matrix polynomial $R_{n,k} (x)$ are defined by
the roots of the polynomial (of degree $np$)
$$
\det R_{n,k} (x),$$ and it can be shown that the orthonormal matrix polynomial $R_{n,k} (x)$ has precisely $np$ real roots, where each root has
at most multiplicity $p$. Throughout this paper let $x_{n,k,j}$ $(j=1, \dots, m)$ denote the different roots of the  matrix orthogonal
polynomial $R_{n,k} (x)$ with multiplicities $\ell_j$, and  consider the empirical distribution of the roots defined by
\begin{eqnarray}\label{2.3}
\delta_{n,k}:=\frac{1}{np}\sum_{j=1}^m\ell_j\delta_{x_{n,k,j}},~n,k\ge 1,
\end{eqnarray}
where $\delta_z$ denotes the Dirac measure at point $z \in \mathbb{R}$. In the following we are interested in the asymptotic properties of this
measure if $n, k \to \infty$. For these investigations we consider sequences $(n_j)_{j \in \mathbb{N}}$, $(k_j)_{j \in \mathbb{N}}$ of positive
integers such that $\frac {n_j}{k_j} \to u$ for some $u \in [0, \infty),$ as $j \to \infty$ and denote the corresponding limit as $\lim_{n/k
\to u}$ (if it exists). Our main result in this section establishes the weak convergence of the sequence of measures $\delta_{n,k}$ under
certain conditions on the matrices $A_{n,k}$  and $B_{n,k}$ if $n/k \to u$ (as $n \to \infty, k \to \infty)$ in the above sense. In the
following sections we will approximate the eigenvalues of the random block matrix $G^{(p)}_n$ by the roots of a specific sequence of matrix
orthogonal polynomials and use this result to derive the asymptotic spectral distribution of the random block matrix $G^{(p)}_n$.

The theory of matrix orthogonal polynomials is substantially richer
than the corresponding theory for the one-dimensional case. Even the case, where all coefficients
in the recurrence relation (\ref{2.1}) are constant, has not been studied in full
detail.
We refer the interested reader Aptekarev and Nikishin (1983), Geronimo (1982),
Dette and Studden (2002) and the references cited in the introduction among many others.
 Before we state our main result regarding the asymptotic
zero distribution of the measure $\delta_{n,k}$ we mention some
facts which are required for the formulation of the following
theorem.
The matrix Chebyshev polynomials of the first kind are defined recursively by
\begin{eqnarray}\nonumber
tT_0^{A,B}(t)&=&\sqrt{2}AT_1^{A,B}(t)+BT_0^{A,B}(t),\\
\label{2.3a}
tT_1^{A,B}(t)&=&AT_2^{A,B}(t)+BT_1^{A,B}(t)+\sqrt{2}AT_0^{A,B}(t),\\
tT_n^{A,B}(t)&=&AT_{n+1}^{A,B}(t)+BT_n^{A,B}(t)+AT_{n-1}^{A,B}(t),~n\ge 2,\nonumber
\end{eqnarray}
where $A$ is a symmetric and non-singular $p \times p$ matrix and $B$ is a symmetric $p \times p$ matrix, and $T_0^{A, B} (t) = I_p$.
Similarly, the Chebyshev polynomials of the second kind are defined by the recursion
\begin{eqnarray}\label{2.3b}
tU_n^{A,B}(t)&=&A^TU_{n+1}^{A,B}(t)+BU_n^{A,B}(t)+AU_{n-1}^{A,B}(t),~n\ge 0,
\end{eqnarray}
with initial conditions $U_{-1}^{A, B} (t) = 0$, $U_0^{A, B} (t) = I_p$. Note that in the case $p=1$, $A=1, B=0$ the polynomials $T^{A,B}_n
(t)$ and $U^{A,B}_n (t)$ are proportional to the classical Chebyshev polynomials of the first and second kind, that is
$$ T^{1,0}_n (t) = \sqrt{2} \cos (n\mbox{arcos} \frac {t}{2}); \ U^{1,0}_n (t) = \frac {\sin \left( (n+1) \arccos \frac {t}{2}\right)}{\sin (\arccos \frac {t}{2})}
\: .$$ In the following discussion the $p \times p$ matrices $A$ and $B$ will depend on a real parameter $u$, i.e.\ $A=A(u), B=B(u)$, and
corresponding matrix measures of orthogonality are denoted by $X_{A(u), B(u)}$ (for the polynomials of the first kind) and $W_{A(u), B(u)}$
(for the polynomials of the second kind). These measures will be normalized such that
$$\int dX_{A(u),B(u)}(t)=\int dW_{A(u),B(u)}(t)=I_p.$$
If $\beta_1(u)\le \ldots \le\beta_p(u)$ are the eigenvalues of the matrix $B(u)$, then it follows from Lemma 2.4 in Duran (1999) that the
matrix
\begin{eqnarray}\label{2.3c1}
K_{A(u),B(u)}(z):=(B(u)-zI_p)^{1/2}A^{-1}(u)(B(u)-zI_p)^{1/2}
\end{eqnarray}
can be diagonalized except for finitely many $z \in \mathbb{C}\backslash [\beta_1(u),\beta_p(u)]$. Throughout this paper the union of the set
of these finitely many complex numbers and the set $[\beta_1, (u), \beta_p (u)]$ will be denoted by $\Delta$ (note that the set $\Delta$
depends on the parameter $u >0$ although this is not reflected in our notation). In this case we have
\begin{eqnarray}\label{2.3c}
K_{A(u),B(u)}(z)=U(z)D(z)U^{-1}(z),~z\in\ce\backslash\Delta,
\end{eqnarray}
 where
$$D(z)=\mbox{diag}\left(\lambda_1^{A(u),B(u)}(z),\ldots,\lambda_p^{A(u),B(u)}(z)\right)
$$
denotes the diagonal matrix of eigenvalues of $K_{A(u), B(u)}(z)$ and $U(z)$ is a unitary $p \times p$ matrix [note that the matrices $U(z)$ and
$D(z)$ depend on the parameter $u$, although this is not reflected by our notation]. It is shown in Duran, Lopez-Rodriguez and Saff (1999) that
under the assumption of a positive definite matrix $A(u)$ and a symmetric matrix $B(u)$ the measure $X_{A(u),B(u)}$ is absolute continuous with
respect to the Lebesgue measure multiplied with the identity matrix and has density
\begin{eqnarray}\label{2.3d}
dX_{A(u),B(u)}(x)&=&A^{-1/2}(u) U(x){\tilde T}(x)U^T(x)A^{-1/2}(u)dx
\end{eqnarray}
where ${\tilde T}(x):=\mbox{diag}({\tilde t}_{11}(x),\ldots,{\tilde t}_{pp}(x))$ denotes a diagonal matrix with elements
\begin{eqnarray}\label{2.3e}
 {\tilde t}_{ii}(x) = \left\{\begin{array}{ll}
       \frac{1}{\pi\sqrt{4-\left( \lambda_i^{A(u),B(u)}(x)\right)^2}}, & \mbox{~if~}{\lambda}_i^{A(u),B(u)}(x)\in (-2,2) \\ \\
             0, & \mbox{~else~}             \end{array}       \right. \qquad i=1,\ldots,p.\end{eqnarray}

For the sake of a simple notation this density is also denoted by $X_{A(u),B(u)} (x)$.

\bigskip

\textbf{Theorem 2.1.} \textit{Consider the sequence of matrix orthonormal polynomials defined by the three term recurrence relation
(\ref{2.1}), where for all $\ell \in \mathbb{N}_0$ and a given $u> 0$ \be
 \label{2.4}
\lim_{\frac{n}{k}\to u}A_{n-\ell,k}=A(u),
\ee
\be
 \label{2.5}
\lim_{\frac{n}{k}\to u}B_{n-\ell,k}=B(u) , \ee with non-singular and symmetric matrices $\{ A (u) | u > 0 \}$ and symmetric matrices $\{ B (u)
| u > 0 \}$. If there exists a number $M > 0$ such that \be \label{2.6} \bigcup^\infty_{k=1} \bigcup^\infty_{n=0} \bigg \{ z \in \mathbb{C} |
\det R_{n,k} (z) = 0  \bigg \} \subset  \ [-M,M] \: , \ee then the empirical measure $\delta_{n,k}$ defined by (\ref{2.3}) converges weakly to
a matrix measure which is absolute continuous with respect to the Lebesgue measure multiplied with the identity matrix. The density of the
limiting distribution is given by \be \label{2.6a}  \frac{1}{u}\int_0^u \mbox{tr}\left[\frac{1}{p}X_{A(s),B(s)}\right]ds ,
 \ee
where $X_{A(s), B(s)}$ denotes the density of the matrix measure corresponding to the matrix Chebyshev polynomials of the first kind.}

\bigskip

{\bf Proof.}  To be precise, let $( R_{n,k} (x) )_{n \ge 0} \ (k \in \mathbb{N})$ denote the sequence of matrix orthonormal  polynomials
defined by the recursive relation (\ref{2.1}) and denote by $x_{n+1,k,1}, \dots, x_{n+1,k,m}$ the different roots of the $(n+1)$th polynomial
$R_{n+1,k} (x)$. At the end of the proof we will show the following auxiliary result, which generalizes a corresponding statement for the case
$p=1$ proved by Kuijlaars and van Assche (1999).

\bigskip

{\bf Lemma A.1.} { \it If all roots  of the matrix orthogonal polynomial $R_{n+1,k} (x)$ defined by (\ref{2.1}) are located in the interval
$[-M,M]$, then the inequality
\begin{eqnarray}\label{2.6aa}
\left|v^T R_{n,k}(z)R^{-1}_{n+1,k}(z) A^{-1}_{n+1,k}v\right|  \le \frac{1}{\mbox{dist}(z,[-M,M])}v^Tv
\end{eqnarray}
 holds for all vectors $v\in\ce^p$, and all complex numbers  $z\in\ce\backslash [-M,M]$. Moreover, we have
 \begin{eqnarray}\label{2.6b}
\left|v^T R_{n,k}(z)R^{-1}_{n+1,k}(z) A^{-1}_{n+1,k}v\right|>\frac{1}{2|z|}v^Tv
\end{eqnarray}
for all vectors $v\in\ce^p$,  $z\in\ce$ with $|z|>M$.}

\bigskip

We now normalize the orthonormal polynomials $R_{n,k} (x)$ such that their leading coefficients are equal to the identity matrix, that is \be
\label{2.7} \underline R_{0,k}(x):=I_p~; \mbox{ ~~ }~\underline R_{n,k}(x):=A_{1,k}\cdots A_{n,k}R_{n,k}(x),~n\ge 1. \ee Then a straightforward
calculation shows that $\underline R_{j,k}^{-1}(x)\underline R_{j+1,k}(x)=R_{j,k}^{-1}(x)A_{j+1,k}R_{j+1,k}(x)$ ($j\ge 0$) and we obtain
\begin{eqnarray}\label{2.8}
\underline R_{n,k}(x)=\prod_{j=0}^{n-1} R_{j,k}^{-1}(x)A_{j+1,k}R_{j+1,k}(x),~n\ge 0.
\end{eqnarray}
This yields
\begin{eqnarray}\label{2.9}
\frac{1}{np} \log \left|\det\underline R_{n,k}(z)\right|&=& 
\frac{1}{np}\sum_{j=0}^{n-1}\log\left|\det\left( R_{j,k}^{-1}(z)A_{j+1,k}R_{j+1,k}(z)\right)\right|\\\nonumber
&=& \frac{1}{p} \int_0^1 \log \left|\det\left( A_{[sn]+1,k}R_{[sn]+1,k}(z)R_{[sn],k}^{-1}(z)\right)\right|ds.
\end{eqnarray}
If $\eta_{n,k,1} (z), \dots, \eta_{n,k,p} (z)$ denote the eigenvalues of the $p \times p$  matrix $R_{n,k}(z)R_{n+1,k}^{-1}(z)A_{n+1,k}^{-1}$,
then it follows
$$
\min_{v\neq 0}\Bigl|\frac{v^TR_{n,k}(z)R_{n+1,k}^{-1}(z)A_{n+1,k}^{-1}v}{v^Tv}  \Bigr |\le \Bigl|\eta_{n,k,i}(z)\Bigr|\le\max_{v\neq
0}\Bigl|\frac{v^TR_{n,k}(z)R_{n+1,k}^{-1}(z)A_{n+1,k}^{-1}v}{v^Tv} \Bigr|
$$
for all $i = 1, \dots, p$ [see Horn and Johnson (1985), p.\ 181]. With these inequalities and Lemma~A.1. we have
\begin{eqnarray}\label{2.10}
\left|\det(R_{n,k}(z)R_{n+1,k}^{-1}(z)A_{n+1,k}^{-1})\right|&=& \Bigl |\prod_{i=1}^p\eta_{n,k,i}(z) \Bigr|
\ge \left(\min_{i=1}^p \left|\eta_{n,k,i}(z)\right|\right)^p\\
\nonumber
&\ge& \left(\min_{v^Tv=1} \left|v^T R_{n,k}(z)R_{n+1,k}^{-1}(z)A_{n+1,k}^{-1}v\right|\right)^p
 >\left(\frac{1}{2|z|}\right)^p \: ,
\end{eqnarray}
whenever $| z| > M $, and
\begin{eqnarray}\label{2.11}
\left|\det(R_{n,k}(z)R_{n+1,k}^{-1}(z)A_{n+1,k}^{-1})\right|
&\le& \left(\max_{i=1}^p
\left|\eta_{n,k,i}(z)\right|\right)^p \\
\nonumber
&\le& \Bigl(\max_{v^Tv=1} \left|v^T R_{n,k}(z)R_{n+1,k}^{-1}(z)A_{n+1,k}^{-1}v\right| \Bigr)^p
 \le\frac{1}{\mbox{dist}(z,[-M,M])^p}.
\end{eqnarray} From (\ref{2.10}) and (\ref{2.11}) we therefore obtain the estimate
\begin{eqnarray}\label{2.12}
\left|\log \left|\det\left(A_{n+1,k}R_{n+1,k}(z)R_{n,k}^{-1}(z)\right)\right|\right|
\le\max\left\{p\left|\log\left(\mbox{dist}(z,[-M,M])\right)\right|,p\left|\log(2|z|)\right|\right\},
\end{eqnarray}
whenever $| z| > M $. Now  the representation (\ref{2.9}) and Lebesgue's theorem yield for $| z | > M$
\begin{eqnarray}\label{2.13}
\lim_{\frac{n}{k}\to u} \frac{1}{np} \log \left|\det\underline R_{n,k}(z)\right|
&=& \frac{1}{p} \int_0^1 \log \bigl|\det( \lim_{\frac{n}{k}\to u}A_{[sn]+1,k}R_{[sn]+1,k}(z)R_{[sn],k}^{-1}(z)) \bigr|ds\\
\nonumber
\nonumber &=& -\frac{1}{p} \int_0^1 \log\left|\det\left(\Phi_{A(su),B(su)}(z)\right)\right|ds,\\\nonumber &=& -\frac{1}{pu} \int_0^u
\log\left|\det\left(\Phi_{A(s),B(s)}(z)\right)\right|ds,
\end{eqnarray}
where
\begin{eqnarray}\label{2.14}
\Phi_{A(u),B(u)}(z):=\int\frac{dW_{A(u),B(u)}(t)}{z-t}
\end{eqnarray}
denotes the (matrix valued) Stieltjes transform of the matrix measure $W_{A(u), B(u)}$ corresponding to the matrix Chebyshev polynomials of the
second kind, and we have used Theorem 2.1 of Duran and Daneri-Vias (2001) for the second equality, which describes the ``ratio asymptotics'' of
the matrices $R_{n+1,k}(z)$ $R^{-1}_{n,k}(z)$. In (\ref{2.13}) the convergence is uniform on compact subsets of $\mathbb{C} \ \backslash \
\Gamma$, where the set $\Gamma$ is defined by
$$
   \Gamma = \bigcap_{m \geq 0} \overline{\bigcup_{j \geq m} \Delta_j} \: ,
   $$
   $\Delta_j = \{ x \in \mathbb{R} \ | \det \ R_{n_j,k_j} (x) = 0 \} $
denotes the set of roots of the matrix polynomial $R_{n_j,k_j} (x),$  and $(n_j, k_j)_{j\ge 0}$ is the sequence along which the convergence is
considered. The Stieltjes transform of the matrix measure $W_{A(u), B(u)}$ has been determined by Duran (1999) as
\begin{eqnarray}\label{2.15}\nonumber
\Phi_{A(u),B(u)}(z)&=&\frac{1}{2}A^{-1}(u)(B(u)-zI_p)^{1/2}\left(-I_p-\sqrt{I_p-4K_{A(u),B(u)}^{-2}(z)}\right)(B(u)-zI_p)^{1/2}A^{-1}(u)\\
\label{stieltjes11}
&=& A^{-1}(u)(B(u)-zI_p)^{1/2}U(z)T(z)U^{-1}(z)(B(u)-zI_p)^{1/2}A^{-1}(u),
\end{eqnarray}
where  the matrix $K_{A(u), B(u)}$ is defined in (\ref{2.3c1}), $ T(z):=\mbox{diag}\left(t_{11}(z),\ldots, t_{pp}(z)\right)$
 is a diagonal matrix with elements given by
$$
  t_{ii}(z) =\frac{-\lambda_i^{A(u),B(u)}(z)-\sqrt{(\lambda_i^{A(u),B(u)}(z))^2-4}}{2\lambda_i^{A(u),B(u)}(z)}$$
$(i=1, \dots, p)$ and $ z \in \mathbb{C} \ \backslash \ \{\mbox{supp}\ (W_{A(u), B(u)}) \cup\Delta\}$. Here and throughout this paper we take
the square root $\sqrt{w}$ such that $| w - \sqrt{w^2-4} | <2$ for $w \in \mathbb{C} \ \backslash \ [-2,2]$, and consequently the function $w -
\sqrt{w^2-4}$ is analytic on $\mathbb{C} \ \backslash \ [-2,2]$. Note that for $z \in \mathbb{C} \ \backslash \ \mbox{supp}\ (W_{A(u), B(u)})$
it follows that $| \lambda_i^{A(u), B(u)} (z) | > 2$. Observing (\ref{2.15}) this implies for the Stieltjes transform
\begin{eqnarray}\label{2.16}
\log|\det(\Phi_{A(u),B(u)}(z))|&=&\log|\det(A^{-2}(u)T(z)(B(u)-zI_p)|\\\nonumber\\\nonumber
&=&\log|\det(A^{-1}(u)T(z)D(z)|
= \log \Bigl |\prod_{j=1}^p\frac{t_{jj}(z)\lambda_j^{A(u),B(u)}(z)}{\alpha_j(u)} \Bigr|\\\nonumber\\\nonumber
&=&\sum_{j=1}^p\log \Bigl |\frac{-\lambda_j^{A(u),B(u)}(z)-\sqrt{(\lambda_j^{A(u),B(u)}
(z))^2-4}}{2\alpha_j(u)} \Bigr| \\
\nonumber
&=& \mbox{Re} \Bigl(\sum_{j=1}^p\log\frac{-\lambda_j^{A(u),B(u)}(z)-\sqrt{(\lambda_j^{A(u),B(u)}(z)
)^2-4}}{2\alpha_j(u)} \Bigr),
\end{eqnarray}
where $\alpha_1 (u), \dots, \alpha_p (u)$ denote the eigenvalues of the matrix $A(u)$. Now define
$$f(z):=\sum_{j=1}^p\log\frac{-\lambda_j^{A(u),B(u)}(z)-\sqrt{(\lambda_j^{A(u),B(u)}(z))^2-4}}{2\alpha_j(u)},$$ then it follows from Kato (1976), p.\
64, that the eigenvalues $\lambda_j^{A(u), B(u)} (z)$ are holomorphic functions on $\mathbb{C} \ \backslash \ \Delta$ and we obtain for any
$$z \in G_0:=\ce\backslash\left\{\mbox{supp}\left(W_{A(u),B(u)}\right)\cup\Delta\right\}$$
that
\begin{eqnarray}\nonumber
\frac{d}{d\bar z}f&=&\sum_{j=1}^p\frac{1}{2f_j(z)\alpha_j(u)} \Bigl( -\frac{d\lambda_j^{A(u),B(u)}(z)}{dx}-\frac{\lambda_j^{A(u),B(u)}
(z)\frac{d\lambda_j^{A(u),B(u)}(z)}{dx}}{\sqrt{(\lambda_j^{A(u),B(u)}(z))^2-4}}
 \Bigr)\\\nonumber\\\nonumber
&+&\sum_{j=1}^pi\frac{1}{2f_j(z)\alpha_j(u)} \Bigl(-\frac{d\lambda_j^{A(u),B(u)}(z)}{dy}-\frac{\lambda_j^{A(u),B(u)}
(z)\frac{d\lambda_j^{A(u),B(u)}(z)}{dy}}{\sqrt{(\lambda_j^{A(u),B(u)}(z))^2-4}}
 \Bigr)\\\nonumber\\\nonumber
&=&\sum_{j=1}^p\frac{1}{2\sqrt{(\lambda_j^{A(u),B(u)}(z))^2-4}}\left(\frac{d}{dx}\lambda_j^{A(u),B(u)}(z)+i\frac{d}{dy}\lambda_j^{A(u),B(u)}(z)\right)=0,
\end{eqnarray}
where the function $f_j$ is defined by
 $$f_j(z)=\frac{-\lambda_j^{A(u),B(u)}(z)-\sqrt{(\lambda_j^{A(u),B(u)}(z)
)^2-4}}{2\alpha_j(u)},~j=1,\ldots,p.$$ This implies that the function $f(z)$ is holomorphic  on $G_0$, and it follows for $z \in G_0$ that
\begin{eqnarray}\label{ablphi}\nonumber
\frac{d}{dz}\log(\det(\Phi_{A(u),B(u)}(z))&=&\frac{1}{2}
\sum_{j=1}^p\frac{1}{\sqrt{(\lambda_j^{A(u),B(u)}(z))^2-4}}\left(\frac{d}{dx}
\lambda_j^{A(u),B(u)}(z)-i\frac{d}{dy}\lambda_j^{A(u),B(u)}(z)\right)\\\nonumber\\
&=&\sum_{j=1}^p\frac{\frac{d}{dz}\lambda_j^{A(u),B(u)}(z)} {\sqrt{(\lambda_j^{A(u),B(u)}(z))^2-4}}.
\end{eqnarray}
In the following discussion define
\begin{eqnarray}\label{2.16a}
U^X(z):=\int\log\frac{1}{|z-t|}\mbox{tr}\left[X_{A(u),B(u)}(t)\right]dt
\end{eqnarray}
as the logarithmic potential of the measure  whose density with respect to the Lebesgue measure is given by $ \mbox{tr} [X_{A(u), B(u)}(t)]$.
Then $U^X (z)$ is harmonic on
$$G_1:=\mathbb{C} \ \backslash \ \mbox{supp} \left(\mbox{tr}[X_{A(u), B(u)}]\right)$$ [see e.g.\ Saff and Totik (1997)].
In the following we show that  the function $f$ satisfies $\mbox{Re }f=U^X+c_1,$ where $c_1 \in \mathbb{C}$ is a constant.
For this purpose we note that the function
\begin{eqnarray}\label{2.16b}
g(z):=\frac{d}{dx}U^X(z)-i\frac{d}{dy}U^X(z)
\end{eqnarray}
is holomorphic on $G_1$ (note that $g$ satisfies the Cauchy-Riemann differential equations  because the logarithmic potential $U^X$ is
harmonic) and satisfies for all $z\in G_1$
\begin{eqnarray}\label{2.16c}
g(z)&=&-\int\frac{\overline{ z-t}}{|z-t|^2}\mbox{tr}\left[X_{A(u),B(u)}(t)\right]dt\\
&=&-\int\frac{\mbox{tr}\left[X_{A(u),B(u)}(t)\right]}{z-t}dt=-\mbox{tr}\left[G_{A(u),B(u)}(z)\right] \: , \nonumber
\end{eqnarray}
where $G_{A(u), B(u)}$ denotes the Stieltjes transform of the matrix measure corresponding to the density $X_{A(u), B(u)} (t)$. In order to
find a representation for the right hand side we note that the function $K_{A(u),B(u)}(z)$ defined in (\ref{2.3c1}) is holomorphic on
$\mathbb{C} \ \backslash \Delta$, and  observing the representation (\ref{2.3c}) we obtain for all $z\in\ce \backslash \Delta$
\begin{eqnarray*} &&-\frac{1}{2}(B(u)-zI_p)^{-1/2}A^{-1}(u)(B(u)-zI_p)^{1/2}-\frac{1}{2}(B(u)-zI_p)^{1/2}A^{-1}(u)(B(u)-zI_p)^{-1/2}\\
&&~=K_{A(u),B(u)}'(z)= U'(z)D(z)U^{-1}(z)+U(z)D'(z)U^{-1}(z)+U(z)D(z)(U^{-1}(z))' \: .
\end{eqnarray*}
>>>>From the identities
\begin{eqnarray*}
A^{-1}(u)(B(u)-zI_p)^{1/2}&=&(B(u)-zI_p)^{-1/2}U(z)D(z)U^{-1}(z),\\
(B(u)-zI_p)^{1/2}A^{-1}(u)&=&U(z)D(z)U^{-1}(z)(B(u)-zI_p)^{-1/2}
\end{eqnarray*}
and
 $(U^{-1}(z))'U(z)=-U^{-1}(z)U'(z)$ it follows that
\begin{eqnarray}\label{2.16d}
&&-\frac{1}{2}U^{-1}(z)(B(u)-zI_p)^{-1}U(z)D(z)-\frac{1}{2}D(z)U^{-1}(z)(B(u)-zI_p)^{-1}U(z)\\\nonumber
&&~=
U^{-1}(z)U'(z)D(z)+D'(z)-D(z)U^{-1}(z)U'(z),
\end{eqnarray}
which yields for the diagonal elements of the matrix $U^{-1}(z)(B(u)-zI_p)^{-1}U(z)$
\begin{eqnarray*}
[U^{-1}(z)(B(u)-zI_p)^{-1}U(z)]_{jj}&=&
\frac{1}{\lambda_j^{A(u),B(u)}(z)}\Bigl(-[U^{-1}(z)U'(z)]_{jj}\lambda_j^{A(u),B(u)}(z)- \frac {d}{dz} \lambda_j^{A(u),B(u)}(z)\\
\nonumber &&~~~~~~+\lambda_j^{A(u),B(u)}(z)[U^{-1}(z)U'(z)]_{jj}\Bigr) \\
\nonumber & =&-\frac{\frac {d}{dz}\lambda_j^{A(u),B(u)}(z)}{\lambda_j^{A(u),B(u)}(z)},
\end{eqnarray*}
for all $z\in\ce\backslash\Delta,~j=1,\ldots,p$. From Zygmunt (2002) and (\ref{2.15}) we have for the Stieltjes transforms $\Phi_{A(u), B(u)}$
and $G_{A(u), B(u)}$ corresponding to the matrix measures $W_{A(u), B(u)}$ and $X_{A(u), B(u)}$ the matrix continued fraction expansion
\begin{eqnarray*}
\nonumber
G_{A(u),B(u)}(z)
&=&\left\{zI_p-B(u)-\sqrt{2}A(u)\left\{zI-B(u)-A(u)\left\{zI_p-B(u)-\ldots\right.\right.\right.\\\nonumber
&&~~~\ldots\left.\left.\left.\right\}^{-1} A(u)\right\}^{-1}\sqrt{2}A(u)\right\}^{-1}\\\label{g1}
&=&\left\{zI_p-B(u)-2A(u)\Phi_{A(u),B(u)}(z)A(u)\right\}^{-1} \\
\label{2.16d1}
&=& (B(u)-zI_p)^{-1/2}U(z)\hat T(z) U^{-1}(z)(B(u)-zI_p)^{-1/2},
\end{eqnarray*}
where $z\in G$,  the set $G$ is given by
$$
G:=\ce\backslash\{\mbox{supp}
\left(W_{A(u),B(u)}\right)\cup\mbox{supp}
\left(X_{A(u),B(u)}\right)\cup \Delta\},
$$
and the diagonal
matrix $\hat T(z)$ is defined by
$${\hat T}(z):= \mbox{diag} (\hat t_{11} (z), \dots, \hat t_{pp} (z))=
\mbox{diag}\Bigl( \frac{\lambda_1^{A(u),B(u)}(z)}{\sqrt{(\lambda_1^{A(u),B(u)}(z))^2-4}}
,\ldots,\frac{\lambda_p^{A(u),B(u)}(z)}{\sqrt{(\lambda_p^{A(u),B(u)}(z))^2-4}}\Bigr ).$$
This yields
\begin{eqnarray}\label{spurg}g(z)
\nonumber
&=& - \mbox{tr}\left[G_{A(u),B(u)}(z)\right]  ~=~ - \sum_{j=1}^p[U^{-1}(z)(B(u)-zI_p)^{-1}U(z)]_{jj}\hat t_{jj}(z)\\\nonumber\\
&=&\sum_{j=1}^p\frac{\frac{d}{dz}\lambda_j^{A(u),B(u)}(z)}{\sqrt{(\lambda_j^{A(u),B(u)}(z))^2-4}},
\end{eqnarray}
for all $z \in G 
.$
Consequently,
we have for all $z\in G$
$$f'(z)=\sum_{j=1}^p \frac{\frac{d}{dz}\lambda_j^{A(u),B(u)}(z)}{\sqrt{(\lambda_j^{A(u),B(u)}(z))^2-4}}=g(z)=\frac{d}{dx}U^X(z)-i\frac{d}{dy}U^X(z).$$
For $h:=\mbox{Re }f$  it follows that
$$
f'=\frac{d}{dx}h-i\frac{d}{dy}h=\frac{d}{dx}U^X-i\frac{d}{dy}U^X,
$$
so that $\frac{d}{dx}(h-U^X)\equiv 0$ and $\frac{d}{dy}(h-U^X)\equiv 0,$ which implies for all $z \in G$ the identity
$\mbox{Re~}f(z)=U^X(z)+c_1$ for some constant $c_1 \in \mathbb{C}$. Therefore it follows that
\begin{equation}
\mbox{Re} \ (f(z)) =
\mbox{Re}\left(\log(\det(\Phi_{A(u),B(u)}(z)))\right)=\log|\det(\Phi_{A(u),B(u)}(z))|
=U^X(z)+c_1
\end{equation}
for any  $z\in G.$
Observing (\ref{2.13}) we finally obtain for all $z \in G$
\begin{eqnarray}\label{2.13d}
\lim_{\frac{n}{k}\to u} \frac{1}{np} \log \left|\det\underline R_{n,k}(z)\right|
&=&-\frac{1}{pu}\int_0^u\int\log\frac{1}{|z-t|}\mbox{tr}\left[X_{A(s),B(s)}(t)\right]dt-c\\\nonumber\\\nonumber
&=&-\int\log\frac{1}{|z-t|}\frac{1}{u}\int_0^u\mbox{tr}\left[\frac{1}{p}X_{A(s),B(s)}(t)\right]dt-c\\\nonumber\\\nonumber &=&-U^\sigma(z)-c,
\end{eqnarray}
where $U^\sigma$ denotes the logarithmic potential of the measure $\sigma$ with Lebesgue density defined in~(\ref{2.6a}) and $c \in \mathbb{C}$
is a constant. Let
$$U^{\delta_{n,k}}(z):=\int\log\frac{1}{|z-t|} \ \delta_{n,k} (dt)$$
denote the logarithmic potential of the measure $\delta_{n,k}$ defined in (\ref{2.3}). From  (\ref{2.13d}) we  obtain for all $z \in G$
\begin{eqnarray}\nonumber\label{2.13e}
\lim_{\frac{n}{k}\to u}U^{\delta_{n,k}}(z)
=\lim_{\frac{n}{k}\to u} \frac {1}{np} \sum^m_{j=1} \log \frac {1}{| z - x_{n,k,j} | ^{\ell_j}}
 &=&\lim_{\frac{n}{k}\to u}\frac{1}{np}\log\frac{1}{|\det\left(R_{n,k}(z)\right)|}\\\nonumber\\
&=&U^\sigma(z)+c.
\end{eqnarray}
The measures in the sequence $(\delta_{n_j,k_j})_{j \in \mathbb{N}}$ have compact support in $[-M,M],$ consequently, $(\delta_{n_j,k_j})_{j \in
\mathbb{N}}$ contains a subsequence which converges weakly to a limit $\mu$ with $\mbox{supp}\ (\mu)\subset [-M,M].$ Therefore we obtain from
(\ref{2.13e})
$$U^\mu(z)=U^\sigma(z)+c,~z\in G,|z|>M,$$ and the assertion of the theorem follows from the fact that the logarithmic  potentials are unique
[see Saff and Totik (1997)].

\bigskip

{\bf Proof of Lemma A.1.} Let $x_{n+1,k,1}, \dots, x_{n+1,k,m}$ denote the different roots of the matrix orthogonal polynomial $R_{n+1,k} (x)$ with
multiplicities $\ell_1, \dots, \ell_m$. Then we obtain from Duran (1996), p.\ 1184, the representation
\begin{eqnarray}\label{a1}
R_{n,k}(z)R^{-1}_{n+1,k}(z) A^{-1}_{n+1,k}=\sum_{j=1}^{m} \frac{C_{n+1,k,j}A_{n+1,k}^{-1}}{z-x_{n+1,k,j}},
\end{eqnarray}
where the weights are given by
\begin{eqnarray}\label{a2}
C_{n+1,k,j}A^{-1}_{n+1,k}&=&R_{n,k}(x_{n+1,k,j})\Gamma_{n+1,k,j}R_{n,k}^T(x_{n+1,k,j}), \\
\label{a3}\nonumber
\Gamma_{n+1,k,j} &=& \frac{\ell_k}{\left(\det (R_{n+1,k}(t))\right)^{(\ell_j)}(x_{n+1,k,j})}\left(\mbox{Adj}(R_{n+1,k}(t)))^{(\ell_j-1)}(x_{n+1,k,j})
Q_{n+1,k}(x_{n+1,k,j}\right),
\end{eqnarray}
and $\mbox{Adj} (A)$ denotes the adjoint of the $p \times p$ matrix $A$, that is
$$
A \ \mbox{Adj} (A) = \mbox{Adj} (A) A = \det (A) I_p\: .
$$
In (\ref{a3}) the matrix polynomial
$$
Q_{n,k} (x) = \int \frac {R_{n,k} (x) - R_{n,k} (t)}{x - t} d \Sigma_k (t)
$$
denotes the first associated matrix orthogonal  polynomial and the matrices $\Gamma_{n,k}$ are nonnegative definite and have rank $\ell_j$ [see
Duran (1996), Theorem 3.1 (2)]. From Duran and Daneri-Vias (2001) we have
\begin{eqnarray}
\label{a5}
\sum_{j=1}^{m} C_{n+1,k,j}A^{-1}_{n+1,k}&=&\sum_{j=1}^m R_{n,k}(x_{n+1,k,j})\Gamma_{n+1,k,j}R_{n,k}^T(x_{n+1,k,j})\\\nonumber &=& \int R_{n,k}(t)\Sigma_k(t)R_{n,k}^T(t)=I_p.
\end{eqnarray}
For any $z\in\ce\backslash [-M,M]$ we obtain the estimate $ |z-x_{n+1,k,j}|\ge \mbox{dist}(z,[-M,M])~\mbox{for all}~j=1,\ldots,m$, because all
roots of the matrix polynomial $R_{n+1,k} (z)$ are located in the interval $[-M,M]$. Note that by the representation (\ref{a2}) the matrix
$C_{n+1,k} A^{-1}_{n+1,k}$ is nonnegative definite and  the representations (\ref{a1}) and (\ref{a5}) yield for any $|z|>M$ and $v \in
\mathbb{C}^p$
\begin{eqnarray*}
\left|v^T R_{n,k}(z)R^{-1}_{n+1,k}(z) A^{-1}_{n+1,k}v\right| &=&\Bigl|v^T
\sum_{j=1}^m\frac{C_{n+1,k,j}A^{-1}_{n+1,k}}{z-x_{n+1,k,j}}v\Bigr|\\\nonumber &\le&
\sum_{j=1}^m\frac{v^TC_{n+1,k,j}A^{-1}_{n+1,k}v}{|z-x_{n+1,k,j}|}\\\nonumber &\le &\frac{1}{\mbox{dist}(z,[-M,M])}\sum_{j=1}^m
v^TC_{n+1,k,j}A^{-1}_{n+1,k}v = \frac{1}{\mbox{dist}(z,[-M,M])}v^Tv,
\end{eqnarray*}
which proves the first inequality of the Lemma. For a proof of the second part we note that for $|z|>M$ we have $|\frac{x_{n+1,k,j}}{z}|<1  $
for all $j=1,\ldots,m$, which gives
$$
\mbox{Re} \Bigl(\frac{1}{1-\frac{x_{n+1,k,j}}{z}}\Bigr)>\frac{1}{2}
$$
for all $j=1,\ldots,m$. With this inequality we obtain
\begin{eqnarray*}
\left|v^T R_{n,k}(z)R^{-1}_{n+1,k}(z) A^{-1}_{n+1,k}v\right|&=&\Bigl|v^T \sum_{j=1}^{m} \frac{C_{n+1,k,j}A_{n+1,k}^{-1}}{z-x_{n+1,k,j}}v\Bigr|\\\nonumber
&=&\frac{1}{|z|} \Bigl|\sum_{j=1}^{m} \frac{v^TC_{n+1,k,j}A_{n+1,k}^{-1}v}{1-\frac{x_{n+1,k,j}}{z}}\Bigr|\\\nonumber
&\ge &\frac{1}{|z|} \mbox{Re}\Bigl(\sum_{j=1}^{m} \frac{v^TC_{n+1,k,j}A_{n+1,k}^{-1}v}{1-\frac{x_{n+1,k,j}}{z}}\Bigr)\\\nonumber
&> &\frac{1}{2|z|}\sum_{j=1}^{m} v^TC_{n+1,k,j}A_{n+1,k}^{-1}v ~=~\frac{1}{2|z|}v^Tv,
\end{eqnarray*}
which proves the second inequality of the Lemma. \hfill $\Box$

\section{Strong and weak asymptotics for eigenvalues of random block matrices}
\def\theequation{3.\arabic{equation}}
\setcounter{equation}{0}

We now consider the random block matrix defined in equation (\ref{1.3}) of the introduction and denote $\tilde\lambda_1^{(n,p)} \leq \dots \leq
\tilde\lambda_n^{(n,p)}$ as its (random) eigenvalues, where $n = mp;~m, p \in \mathbb{N}$. Similarly, define $\tilde x_1^{(n,p)} \leq \dots
\leq \tilde x_n^{(n,p)}$ as (deterministic) eigenvalues of the matrix orthonormal  polynomial $\tilde R_{m,n}^{(p)} (x) \in \mathbb{R}^{p
\times p}$, which is defined recursively by $\tilde R^{(p)}_{0,n} (x) = I_p, \ \tilde R^{(p)}_{-1,n} (x) = 0$, \be \label{3.1} x \ \tilde
R_{m,n}^{(p)} (x) = \tilde A_{m+1,n}^{(p)} \ \tilde R_{m+1,n}^{(p)} (x) + \tilde B_{m,n}^{(p)} \ \tilde R_{m,n}^{(p)} (x) + \tilde
A_{m,n}^{(p)^T} \ \tilde R_{m-1,n}^{(p)} (x) \ee $(m \geq 0)$, where the $p \times p$ block matrices $\tilde A_{i,n}^{(p)}$ and $\tilde
B_{i,n}^{(p)}$ are given by {\scriptsize\begin{eqnarray}\label{3.2}
{\tilde A}_{i,n}^{(p)}=\frac{1}{\sqrt{2}}\left(
\begin{array}{cccccccc}
\sqrt{((i-1)p+1)\gamma_p} &  \sqrt{((i-1)p+2)\gamma_{p-1}}  & \sqrt{((i-1)p+3)\gamma_{p-2}}  & \cdots   &  \sqrt{ip\gamma_{1}}       \\
\sqrt{((i-1)p+2)\gamma_{p-1}} &\sqrt{((i-1)p+2)\gamma_p}&  \sqrt{((i-1)p+3)\gamma_{p-1}}  & \cdots &    \sqrt{ip\gamma_{2}}  \\
  &  &  &  &   \\
 \vdots &\ddots & \ddots & \ddots & \vdots  \\
  & & &  &   \\
 \sqrt{(ip-1)\gamma_2}  &\cdots  &  \sqrt{(ip-1)\gamma_{p-1} }   & \sqrt{(ip-1)\gamma_p} & \sqrt{ip\gamma_{p-1}} \\
 \sqrt{ip\gamma_1}  & \cdots&  \sqrt{ip\gamma_{p-2}}&\sqrt{ip\gamma_{p-1}} & \sqrt{ip\gamma_p} \\
 \end{array}
\right),
\end{eqnarray}}
$i\ge 1,$ and
{\scriptsize \begin{eqnarray}\label{3.3}
{\tilde B}_{i,n}^{(p)}=\frac{1}{\sqrt{2}}\left(
\begin{array}{cccccccc}
0 &  \sqrt{(ip+1)\gamma_1} &\sqrt{(ip+1)\gamma_2}  & \cdots   &   \sqrt{(ip+1)\gamma_{p-1}}       \\
\sqrt{(ip+1)\gamma_1}  &0& \sqrt{(ip+2)\gamma_1} & \cdots &    \sqrt{(ip+2)\gamma_{p-2}}  \\
  &  &  &  &   \\
 \vdots &\ddots & \ddots & \ddots & \vdots  \\
  & & &  &   \\
 \sqrt{(ip+1)\gamma_{p-2}}  &\cdots  &  \sqrt{((i+1)p-2)\gamma_{1}}   &  0 & \sqrt{((i+1)p-2)\gamma_1} \\
 \sqrt{(ip+1)\gamma_{p-1} }  & \cdots&  \sqrt{((i+1)p-2)\gamma_{2}} & \sqrt{((i+1)p-1)\gamma_1} & 0 \\
 \end{array}
\right),
\end{eqnarray}}
$i\ge 0,$ respectively, and $\gamma_1,\ldots,\gamma_p >0$ are given constants. Our first result provides a strong uniform approximation of the
ordered random eigenvalues of the matrices $G^{(p)}_n$ by the ordered deterministic roots of the matrix orthogonal polynomials $R^{(p)}_{m,n}
(x)$.

\bigskip

{\bf Theorem 3.1.} {\it Let $\tilde \lambda_{1}^{(n,p)} \leq \dots \leq \tilde \lambda_n^{(n,p)}$ denote the
 eigenvalues of the random matrix $G_n^{(p)}$ defined by (\ref{1.3}) and $\tilde x_{1}^{(n,p)} \leq \dots \leq \tilde x_n^{(n,p)}$
denote the roots of the matrix orthonormal polynomial $\tilde R_{m,n}^{(p)} (x)$ defined by  the recurrence relation (\ref{3.1}), then
\begin{eqnarray}\label{3.5}
\max_{1\le j\le n}\left|\tilde \lambda_j^{(n,p)}-\tilde  x_j^{(n,p)}\right|\le\left[{\log n}\right]^{1/2}S~~\forall~n\ge 2,
\end{eqnarray}
where $S$ is a random variable such that $S < \infty  \ a.s.$ }

\bigskip

{\bf Proof.} For the proof we establish the bound
\begin{eqnarray}\label{3.6}
P\left\{\max_{1\le j\le n}\left| \tilde \lambda_j^{(n,p)}- \tilde x_j^{(n,p)}\right|\ge \eps\right\}\le 2n(p+1) \exp\left(\frac{-\eps^2}{18p^2}\right),
\end{eqnarray}
then the assertion follows along the lines of the proof of Theorem 2.2 in Dette and Imhof (2007). First we note that the roots of the $m$th
matrix orthonormal polynomial $\tilde R^{(p)}_{m,n} (x)$ are the eigenvalues of the tridiagonal block matrix
\begin{eqnarray}\label{3.6a}
{\tilde F}_n^{(p)}:=\left(
\begin{array}{cccccccc}
{\tilde B}_{0,n} &  {\tilde A}_{1,n} &  &    &          \\
{\tilde A}_{1,n} & {\tilde B}_{1,n} &  {\tilde A}_{2,n} &   &   \\
 &\ddots  & \ddots & \ddots & &  \\
  &  &  {\tilde A}_{\frac{n}{p}-2,n}  &  {\tilde B}_{\frac{n}{p}-2,n} &  {\tilde A}_{\frac{n}{p}-1,n}  \\
 & & & {\tilde A}_{\frac{n}{p}-1,n} & {\tilde B}_{\frac{n}{p}-1,n}  \\
 \end{array}
\right)\in\er^{n\times n} \: ,
 \end{eqnarray}
where the blocks $\tilde A_{i,n}^{(p)},$ $i=1,\ldots,n/p-1,$ and  $\tilde B_{i,n}^{(p)},$ $i=0,\ldots,n/p-1,$ are defined by (\ref{3.2}) and
(\ref{3.3}), respectively. Moreover, by interchanging rows and columns of the matrix $\tilde F_{n}^{(p)}$ it is easy to see that the matrix
\begin{eqnarray}\label{3.7}
F_n^{(p)}=\left(
\begin{array}{cccccccc}
E_{0,n}^{(p)} &  D_{1,n}^{(p)} &  &    &          \\
D_{1,n}^{(p)} & E_{1,n}^{(p)} &  D_{2,n}^{(p)} &   &   \\
 &\ddots  & \ddots & \ddots & &  \\
  &  &  D_{\frac{n}{p}-2,n}^{(p)}  &  E_{\frac{n}{p}-2,n}^{(p)} &  D_{\frac{n}{p}-1,n}^{(p)}  \\
 & & & D_{\frac{n}{p}-1,n}^{(p)} & E_{\frac{n}{p}-1,n}^{(p)} \\
 \end{array}
\right)
 \end{eqnarray}
 with symmetric blocks

{\scriptsize
 \begin{eqnarray}\label{3.7a}\nonumber
 E_{i,n}^{(p)}=\frac{1}{\sqrt{2}}\left(
\begin{array}{cccccccc}
0 &  \sqrt{\gamma_1 (n-ip-1)} &\sqrt{\gamma_2 (n-ip-2)}  & \cdots   &   \sqrt{\gamma_{p-1} (n-(i+1)p+1)}       \\
\sqrt{\gamma_1 (n-ip-1)}  &0& \sqrt{\gamma_1 (n-ip-2)} & \cdots &    \sqrt{\gamma_{p-2} (n-(i+1)p+1)}  \\
  &  &  &  &   \\
 \vdots &\ddots & \ddots & \ddots & \vdots  \\
  & & &  &   \\
 \sqrt{\gamma_{p-2} (n-(i+1)p+2)}  &\cdots  &  \sqrt{\gamma_{1} (n-(i+1)p+1)}   &  0 & \sqrt{\gamma_1 (n-(i+1)p+1)} \\
 \sqrt{\gamma_{p-1} (n-(i+1)p+1)}  & \cdots&  \sqrt{\gamma_{2} (n-(i+1)p+1)} & \sqrt{\gamma_1 (n-(i+1)p+1)} & 0 \\
 \end{array}
\right),
 \end{eqnarray}}
$i=0,\ldots,\frac{n}{p}-1,$ and

{\scriptsize
 \begin{eqnarray}\label{3.7b}\nonumber
D_{i,n}^{(p)}=\frac{1}{\sqrt{2}}\left(
\begin{array}{cccccccc}
\sqrt{\gamma_p(n-ip)} &  \sqrt{\gamma_{p-1}(n-ip)}  & \sqrt{\gamma_{p-2}(n-ip)}  & \cdots   &  \sqrt{\gamma_{1}(n-ip)}       \\
\sqrt{\gamma_{p-1}(n-ip)} &\sqrt{\gamma_p(n-ip-1)}&  \sqrt{\gamma_{p-1}(n-ip-1)}  & \cdots &    \sqrt{\gamma_{2}(n-ip-1)}  \\
  &  &  &  &   \\
 \vdots &\ddots & \ddots & \ddots & \vdots  \\
  & & &  &   \\
 \sqrt{\gamma_2(n-ip)}  &\cdots  &  \sqrt{\gamma_{p-1} (n-(i+1)p+3)}   & \sqrt{\gamma_p (n-(i+1)p+2)} & \sqrt{\gamma_{p-1} (n-(i+1)p+2)} \\
 \sqrt{\gamma_1(n-ip)}  & \cdots&  \sqrt{\gamma_{p-2} (n-(i+1)p+3)}&\sqrt{\gamma_{p-1} (n-(i+1)p+2)} & \sqrt{\gamma_p (n-(i+1)p+1)} \\
 \end{array}
\right),
 \end{eqnarray}}
$(i=1,\ldots,\frac{n}{p}-1)$ also has eigenvalues $\tilde x_{1}^{(n,p)} \leq \dots \leq \tilde x_n^{(n,p)}$. From Horn and Johnson (1985) we
therefore obtain the estimate
\begin{eqnarray}\label{3.8}
\max_{1\le j\le n} \left| \tilde
\lambda_j^{(n,p)}-\tilde x_j^{(n,p)}\right|
\le \left\|G_n^{(p)}-F_n^{(p)}\right\|_\infty,
\end{eqnarray}
where for a matrix $A = (a_{ij})_{i,j=1,\dots,n}$
$$
\left\|A\right\|_\infty=\max_{1\le i\le n}\sum_{j=1}^n\left|a_{ij}\right |$$
denotes the maximum row sum norm. If
$$
Z_n:=\max\left\{\max_{1\le j\le n} \left|N_j\right|,\max_{1\le j\le n-1} \frac{1}{\sqrt{2}}\left|X_{j\gamma_1}-\sqrt{j\gamma_1}\right|,
\ldots,\max_{1\le j\le n-p} \frac{1}{\sqrt{2}}\left|X_{j\gamma_p}-\sqrt{j\gamma_p}\right|\right\},
$$
then it follows for $i=1, \dots, p, n-p+1, \dots, n$ by a straightforward but tedious calculation
$$
\sum_{j=1}^n\left|\{G_n^{(p)}-F_n^{(p)}\}_{ij}\right|\le 2pZ_n,
$$
and for $i = p+1, \dots, n-p$
$$
\sum_{j=1}^n\left|\{G_n^{(p)}-F_n^{(p)}\}_{ij}\right|\le 3pZ_n,
$$
which yields observing (\ref{3.8})
$$
\max_{1\le j\le n} \left| \tilde \lambda_j^{(n,p)}-\tilde x_j^{(n,p)}\right|\le 3pZ_n
$$
and
\begin{eqnarray}\label{3.9}
P\left\{\max_{1\le j\le n}\left| \tilde  \lambda_j^{(n,p)}- \tilde  x_j^{(n,p)}\right|\ge\eps\right\}\le P\left\{Z_n\ge \frac{\eps}{3p}\right\}.
\end{eqnarray} From Dette and Imhof (2007) we  obtain the estimates
$$P\Bigl\{\frac{\left|X_{j\gamma_k}-\sqrt{j\gamma_k}\right|}{\sqrt{2}}\ge\frac{\eps}{3p}\Bigr\}\le 2e^{-\eps^2/9p^2}~,~
P\Bigl\{\max_{1\le j\le n} |N_j|\ge \frac{\eps}{3p}\Bigr\} \le 2n e^{-\eps^2/18p^2} \: ,
$$
and a combination of these inequalities with (\ref{3.9}) yields
$$
P\left\{\max_{1\le j\le n}|\tilde \lambda_j^{(n,p)}-\tilde x_j^{(n,p)}|\ge
\eps\right\}\le P\left\{Z_n\ge \frac{\eps}{3p}\right\} \le
2n(p+1)e^{-\eps^2/18p^2},
$$
which completes the proof of the theorem. \hfill $\Box$

\bigskip

{\bf Theorem 3.2.} {\it Let $\lambda_1^{(n,p)} \leq \dots \leq \lambda_n^{(n,p)}$ denote the eigenvalues of the random matrix $\frac
{1}{\sqrt{n}} G_n^{(p)}$, where $G_n^{(p)}$ is defined in (\ref{1.3}) and $\gamma_1, \dots, \gamma_p > 0$ are chosen such that all blocks
$D_i^{(p)}$ in the matrix $F^{(p)}_n$ defined by (\ref{3.7}) are non-singular. If \be \sigma_n = \frac {1}{n} \sum^n_{j=1}
\delta_{\lambda_j^{(n,p)}} \ee denotes the empirical distribution of the eigenvalues of the matrix $\frac {1}{\sqrt{n}} G_n^{(p)}$, then
$\sigma_n$ converges weakly to a measure, which is absolute continuous with respect to the Lebesgue measure. The density of this measure is
given by \be \label{3.9a} f(t) = \int_0^\frac{1}{p}\mbox{tr}\left[X_{A^{(p)}(s),B^{(p)}(s)}(t)\right]ds, \ee where $X_{A^{(p)}(s),
B^{(p)}(s)}(t)$ denotes the Lebesgue density of the matrix measure  corresponding to the Chebyshev polynomials of the first kind defined in
(\ref{2.3a}) with matrices
\begin{eqnarray}\nonumber
A^{(p)}(s)&:=&\sqrt{\frac{sp}{2}}\left(
\begin{array}{cccccccc}
\sqrt{\gamma_p} &  \sqrt{\gamma_{p-1}}  & \sqrt{\gamma_{p-2}}  & \cdots   &  \sqrt{\gamma_{1}}       \\
\sqrt{\gamma_{p-1}} &\sqrt{\gamma_p}&  \sqrt{\gamma_{p-1}}  & \cdots &    \sqrt{\gamma_{2}}  \\
  &  &  &  &   \\
 \vdots &\ddots & \ddots & \ddots & \vdots  \\
  & & &  &   \\
 \sqrt{\gamma_2}  &\cdots  &  \sqrt{\gamma_{p-1} }   & \sqrt{\gamma_p} & \sqrt{\gamma_{p-1}} \\
 \sqrt{\gamma_1}  & \cdots&  \sqrt{\gamma_{p-2}}&\sqrt{\gamma_{p-1}} & \sqrt{\gamma_p} \\
 \end{array}
\right)\in\er^{p\times p}, \\
\nonumber  \\
B^{(p)}(s)&:=&\sqrt{\frac{sp}{2}}\left(
\begin{array}{cccccccc}
0 &  \sqrt{\gamma_1} &\sqrt{\gamma_2}  & \cdots   &   \sqrt{\gamma_{p-1}}       \\
\sqrt{\gamma_1}  &0& \sqrt{\gamma_1} & \cdots &    \sqrt{\gamma_{p-2}}  \\
  &  &  &  &   \\
 \vdots &\ddots & \ddots & \ddots & \vdots  \\
  & & &  &   \\
 \sqrt{\gamma_{p-2}}  &\cdots  &  \sqrt{\gamma_{1}}   &  0 & \sqrt{\gamma_1} \\
 \sqrt{\gamma_{p-1} }  & \cdots&  \sqrt{\gamma_{2}} & \sqrt{\gamma_1} & 0 \\
 \end{array}
\right)\in\er^{p\times p}.
 \end{eqnarray}}

\bigskip

{\bf Proof.} Let $R^{(p)}_{m,n} (x)$ $(m=n/p)$ denote the orthonormal polynomials satisfying the recurrence relation (\ref{2.1}) with coefficients $A^{(p)}_{i,n} = \frac {1}{\sqrt{n}}  \ \tilde A^{(p)}_{i,n}$ and $B^{(p)}_{i,n} = \frac {1}{\sqrt{n}}  \ \tilde B^{(p)}_{i,n}$,
where $\tilde A^{(p)}_{i,n}$ and $\tilde B^{(p)}_{i,n}$ are defined by (\ref{3.2}) and (\ref{3.3}), respectively, then we have
\begin{eqnarray*}
R^{(p)}_{m,n} (x) &=& \tilde R^{(p)}_{m,n} (\sqrt{n} x), \\
x_j^{(n,p)} &=& \frac {\tilde x_j^{(n,p)}}{\sqrt{n}}, \quad \quad  j=1, \dots, n,
\end{eqnarray*}
where the matrix  polynomial $\tilde R^{(p)}_{m,n}(x)$ is defined by (\ref{3.1}) with corresponding roots $\tilde x^{(n,p)}_j$,  and
$x^{(n,p)}_j$ denotes the $j$th root of the matrix polynomial $R^{(p)}_{m,n} (x)$. From the definition of the matrices $\tilde A^{(p)}_{i,n}$
and $\tilde B^{(p)}_{i,n}$ in (\ref{3.2}) and (\ref{3.3}) we have
\begin{eqnarray}\nonumber
\lim_{n\to\infty}{ A}_{\frac{n}{p}-1,n}^{(p)}&=&\frac{1}{\sqrt{2}}\left(
\begin{array}{cccccccc}
\sqrt{\gamma_p} &  \sqrt{\gamma_{p-1}}  & \sqrt{\gamma_{p-2}}  & \cdots   &  \sqrt{\gamma_{1}}       \\
\sqrt{\gamma_{p-1}} &\sqrt{\gamma_p}&  \sqrt{\gamma_{p-1}}  & \cdots &    \sqrt{\gamma_{2}}  \\
  &  &  &  &   \\
 \vdots &\ddots & \ddots & \ddots & \vdots  \\
  & & &  &   \\
 \sqrt{\gamma_2}  &\cdots  &  \sqrt{\gamma_{p-1} }   & \sqrt{\gamma_p} & \sqrt{\gamma_{p-1}} \\
 \sqrt{\gamma_1}  & \cdots&  \sqrt{\gamma_{p-2}}&\sqrt{\gamma_{p-1}} & \sqrt{\gamma_p} \\
 \end{array}
\right)=:A^{(p)}, \\
\nonumber  \\
\nonumber
\lim_{n\to\infty} 
{ B}_{\frac{n}{p}-1,n}^{(p)}&=&\frac{1}{\sqrt{2}}\left(
\begin{array}{cccccccc}
0 &  \sqrt{\gamma_1} &\sqrt{\gamma_2}  & \cdots   &   \sqrt{\gamma_{p-1}}       \\
\sqrt{\gamma_1}  &0& \sqrt{\gamma_1} & \cdots &    \sqrt{\gamma_{p-2}}  \\
  &  &  &  &   \\
 \vdots &\ddots & \ddots & \ddots & \vdots  \\
  & & &  &   \\
 \sqrt{\gamma_{p-2}}  &\cdots  &  \sqrt{\gamma_{1}}   &  0 & \sqrt{\gamma_1} \\
 \sqrt{\gamma_{p-1} }  & \cdots&  \sqrt{\gamma_{2}} & \sqrt{\gamma_1} & 0 \\
 \end{array}
\right)=:B^{(p)} \: ,
 \end{eqnarray}
and Gerschgorin's disc theorem implies that all roots of the polynomials $R^{(p)}_{m,n} (x)$ are located in a compact interval, say $[-M,M]$
[see also the proof of Lemma 2.1 in Duran (1999)]. Moreover, we have for any $\ell \in \mathbb{N}_0$
\begin{eqnarray*}
\lim_{\frac{i}{n}\to u} {B}_{i-\ell,n}^{(p)}=\sqrt{up}
B^{(p)}=:B^{(p)}(u) , \\
\lim_{\frac{i}{n}\to u} {A}_{i-\ell,n}^{(p)}=\sqrt{up} A^{(p)}=:A^{(p)}(u),
 \end{eqnarray*}
where $u > 0$ and the matrix $A^{(p)} (u)$ is non-singular by assumption. Consequently, Theorem~2.1 is applicable and yields (note that
$\lim_{n \to \infty} (n/p)/n = 1/p)$ that the empirical distribution of the roots of the matrix polynomials $R^{(p)}_{m,n} (x) \ (m=n/p)$ \be
\delta_n = \frac {1}{n} \sum^n_{j=1} \delta_{x_j^{(n,p)}} \ee converges weakly to the measure with Lebesgue density $f$ defined in
(\ref{3.9a}). Next we use Theorem~3.1 which shows that \be \max^n_{j=1} \ | \ \lambda_j^{(n,p)} - x_j^{(n,p)} \ | \ = \frac {1}{\sqrt{n}}
\max^n_{j=1} \ | \ \tilde\lambda_j^{(n,p)} - \tilde x_j^{(n,p)} \ | \ = O \Bigl( \Bigl( \frac {\log n}{n} \Bigr)^{1/2}\Bigr) \ee almost surely.
Therefore we obtain for the Levy distance $L$ between the distribution functions $F_{\sigma_n}$ and $F_{\delta_n}$ of the measures $\sigma_n$
and $\delta_n$ \be L^3 (F_{\sigma_n}, F_{\delta_n}) \leq \frac {1}{n} \sum^n_{j=1} \ | \ \lambda_j^{(n,p)} - x_j^{(n,p)} \ |^2 = O \left( \frac
{\log n}{n}\right) \ee almost surely [for the inequality see Bai (1999), p.\ 615]. Consequently, it follows that the spectral measure
$\sigma_n$ of the matrix $\frac {1}{\sqrt{n}} G^{(p)}_n$ converges also weakly with the same limit as $\delta_n$, that is the measure with
Lebesgue density $f$ defined by (\ref{3.9a}).

 \hfill $\Box$

 \section{Examples}
\def\theequation{4.\arabic{equation}}
\setcounter{equation}{0}

We conclude this paper with a discussion of a few examples. First note that in the case $p=1$ Theorem 3.2 yields for the limiting distribution
of eigenvalues of the matrix $\frac {1}{\sqrt{n}}G_n^{(1)}$ the Wigner's semicircle law with density
$$
f(x)=\frac{1}{\pi\gamma_1}\sqrt{2\gamma_1-x^2}I_{\left\{\sqrt{-2\gamma_1}< x< \sqrt{2\gamma_1}\right\}},
$$
which have been
considered by numerous authors.

Next we concentrate on the case $p=2$, for which it is easily seen that the matrix $D_i^{(2)}$ in (\ref{3.7a}) is non-singular whenever
$\gamma_2 \neq \gamma_1$. In this case the density of the limit of the spectral measure is a mixture of two arcsine densities and given by
\begin{eqnarray}\nonumber
f(x) =
\lim_{n\to\infty}\frac{1}{n}\sum_{j=1}^n\delta_{{\lambda}_j^{(n,2)}}&=&\sum_{j=1}^2\int_0^{1/2}\frac{1}{\pi\sqrt{4\alpha_j^2(s)-\left(x-\beta_j(s)\right)^2}}I_{\left\{-2\alpha_j(s)+\beta_j(s)<
x< 2\alpha_j(s)+\beta_j(s)\right\}}ds
\end{eqnarray}
where
$$
\alpha_{1}(s)=\sqrt{s}\left(\sqrt{\gamma_2}+\sqrt{\gamma_1}\right),~\alpha_{2}(s)
 = \sqrt{s}\left(\sqrt{\gamma_2}-\sqrt{\gamma_1}\right),~\beta_{1}(s)=\sqrt{s\gamma_1},~ \beta_{2}(s)=-\sqrt{s\gamma_1} \: .
$$
In Figure 1 and 2 we display the limiting spectral density corresponding to the case $\gamma_1 = 2, \gamma_2 = 8$
 and $\gamma_1 = 1, \gamma_2=100,$
respectively. The left part of the figures shows a simulated histogram of the eigenvalues of the matrix $\frac {1}{\sqrt{n}} \ G^{(2)}_n$ for
$n = 5000$, while the right part of the figures shows the corresponding limiting distribution obtained from Theorem 3.2.

\begin{figure} \label{fig4.1}
\begin{center}
\end{center}
\caption{ \it Simulated and limiting spectral density of the random block matrix $G_n^{(p)}/\sqrt{n}$ in the case $p=2$, $\gamma_1=2$,
$\gamma_2=8$. In the simulation the eigenvalue distribution of a $5000 \times 5000$ matrix was calculated (i.e. m = n/p = 2500).}
\end{figure}

\begin{figure} \label{fig4.2}
\begin{center}
\end{center}
\caption{ \it Simulated and limiting spectral density of the random block matrix $G_n^{(p)}/\sqrt{n}$ in the case $p=2$, $\gamma_1=1$,
$\gamma_2=100$. In the simulation the eigenvalue distribution of a $5000 \times 5000$ matrix was calculated (i.e. m = n/p = 2500).}
\end{figure}
If $p \geq 3$ the general formulas for the density of the limit distribution are too complicated to be displayed here, but it can be shown that
$$\mbox{tr}[X_{A^{(p)}(u),B^{(p)}(u)}(t)]=\sum_{j=1}^p\frac{-\frac{d}{dt}\lambda_j^{A^{(p)}(u),B^{(p)}(u)}(t)}{\pi\sqrt{4-(\lambda_j^{A^{(p)}(u),B^{(p)}(u)}(t))^2}}I_{\{-2<\lambda_j^{A^{(p)}(u),B^{(p)}{(p)}}(t)<2\}},$$
if the matrix $A^{(p)}(u)$ is positive definite. This identity follows in a similiar way as (\ref{spurg}).
In Figure~3, 4 and 5 we show a simulated histogram of the eigenvalues of $G_n^{(3)}/\sqrt{n}$ and the corresponding density of the limit
distribution obtained from Theorem 3.2 in the case $\gamma_1=\gamma_2=4,~\gamma_3=100;$ $\gamma_1=1,~\gamma_2=4,\gamma_3=25$ and $\gamma_1=1,~\gamma_2=100,\gamma_3=200,$ respectively.

\begin{figure} \label{fig4.3}
\begin{center}
\end{center}
\caption{ \it Simulated and limiting spectral density of the random block matrix $G_n^{(p)}/\sqrt{n}$ in the case $p=3$, $\gamma_1=\gamma_2=4$,
$\gamma_3=100$. In the simulation the eigenvalue distribution of a $5001 \times 5001$ matrix was calculated. (i.e. m = n/p = 1667)}
\end{figure}

\begin{figure} \label{fig4.4}
\begin{center}
\end{center}
\caption{ \it Simulated and limiting spectral density of the random block matrix $G_n^{(p)}/\sqrt{n}$ in the case $p=3$,
$\gamma_1=1,~\gamma_2=4$, $\gamma_3=25$. In the simulation the eigenvalue distribution of a $5001 \times 5001$ matrix was calculated (i.e. m =
n/p = 1667).}
\end{figure}

\begin{figure} \label{fig4.5}
\begin{center}
\end{center}
\caption{ \it Simulated and limiting spectral density of the random block matrix $G_n^{(p)}/\sqrt{n}$ in the case $p=3$,
$\gamma_1=1,~\gamma_2=100$, $\gamma_3=200$. In the simulation the eigenvalue distribution of a $5001 \times 5001$ matrix was calculated (i.e. m
= n/p = 1667).}
\end{figure}

\vskip 1cm

{\bf Acknowledgements.} The authors are grateful to Martina Stein  who typed parts of this paper with
considerable technical expertise.
The work of the authors was supported by the Sonderforschungsbereich Tr/12,
Fluctations and universality of invariant random matrix ensembles (project C2) and
in part by a NIH grant award
IR01GM072876:01A1.

\newpage

\vskip 1cm

\section*{References}
\smallskip
\bigskip

A.I.\ Aptekarev, E.M.\ Nikishin (1983). The scattering problem for discrete Sturm-Liouville operator.
Mat.\ Sb.\ 121 (163), 327-358; Math.\ USSR\ Sb.\ 49, 325-355, 1984.

\smallskip

Z.D.\ Bai (1999). Methodologies in spectral analysis of large-dimensional random matrices,
a review. Statist. Sinica 9, 611-677.

\smallskip

H.\ Dette and L.A.\ Imhof (2007). Uniform approximation of eigenvalues in Laguerre and Hermite $ \beta$-ensembles by
roots of orthogonal polynomials. Trans. Amer.\ Math.\ Soc.\ 359, 4999-5018.

\smallskip

H.\ Dette and W.J.\ Studden (2002). Matrix measures, moment spaces and Favard's theorem on the interval $[0,1]$
and $[0,\infty]$. Lin.\ Alg.\ Appl.\ 345, 169-193.

\smallskip

I.\ Dumitriu and A.\ Edelman (2002). Matrix models for beta ensembles. J.\ Math.\ Phys.\ 43, 5830-5847.

\smallskip

I.\ Dumitriu and A.\ Edelman (2005). Eigenvalues of Hermite and Laguerre ensembles: large beta asymptotics.
Ann.\ Inst. Poincar\'{e}, Probab.\ Statist. 41, 1083-1099.
\smallskip

A.J.\ Duran (1996). Markov's theorem for orthogonal matrix polynomials, Can.\ J.\ Math.\ 48, 1180-1195.

\smallskip

A.J.\ Duran (1999). Ratio asymptotics for orthogonal matrix polynomials,
J.\ Approx.\ Theory \ 100, 304-344.

\smallskip

A.J.\ Duran and E.\ Daneri-Vias (2001). Ratio asymptotics for orthogonal matrix polynomials with unbounded
recurrence coefficients, J.\ Approx.\ Theory 110, 1-17

\smallskip

A.J.\ Duran and P.\ Lopez-Rodriguez (1996). Orthogonal matrix polynomials: zeros and Blumenthal's
theorem, J.\ Approx.\ Theory \ 84, 96-118.

\smallskip

A.J.\ Duran, P.\ Lopez-Rodriguez and E.B.\ Saff (1999). Zero asymptotic behaviour for orthogonal matrix polynomials, J.\ Anal.\ Math. 78, 37-60.

\smallskip

F.J.\ Dyson (1962). The threefold way. Algebraic structure of symmetry groups and ensembles in
quantum mechanics. J.\ Math.\ Phys.\ 3, 1199-1215.

\smallskip

J.S.\ Geronimo (1982). Scattering theory and matrix orthogonal polynomials on the real line.
Circuits Systems Signal Process 1, 471-495.

\smallskip

V.L.\ Girko (2000). Random block matrix density and SS-law. Random Oper. Stochastic Equations\ 8, 189-194.

\smallskip

R.A.\ Horn and C.R.\ Johnson (1985). Matrix Analysis, Cambridge University Press, Cambridge.

\smallskip

T.\ Kato (1976). Perturbation Theory for Linear Operators. Springer, Berlin.

\smallskip

A.B.J.\ Kuijlaars and W.\ van Assche (1999). The asymptotic zero distribution of orthogonal polynomials
with varying recurrence coefficients. J.\ Approx.\ Theory 99, 167-197.

\smallskip

M.L.\ Mehta (1967). Random Matrices. Academic Press, New York.

\smallskip

T.\ Oraby (2006a). The limiting spectra of Girko's block-matrix.
Front for the arXiv, math.PR Probability Theory, arXiv:math/0612177.

\smallskip

T.\ Oraby (2006b). The spectral laws of hermitian block-matrices with large random blocks.
Front for the arXiv, math.PR Probability Theory, arXiv:0704.2904.

\smallskip

E.B.\ Saff and V.\ Totik (1997). Logarithmic Potentials with External Fields, Springer, Berlin.

\smallskip

A.\ Sinap and W.\ van Assche (1996). Orthogonal matrix polynomials and applications. J.\ Comput.\ Appl.\ Math.\ 66, 27-52.

\smallskip

M.J. Zygmunt (2002). Matrix Chebyshev polynomials and continued fractions. Lin.\ Alg.\ Appl.\ 340, 155-168.

\end{document}